\let\clineorig\cline  % bug fix: see https://tex.stackexchange.com/questions/678195/latex-table-missing-border-lines

\documentclass{eptcs}
\providecommand{\event}{VSS 2025} % Name of the event you are submitting to

\ifpdf
  \usepackage{underscore}         % Only needed if you use pdflatex.
  \usepackage[T1]{fontenc}        % Recommended with pdflatex
\else
  \usepackage{breakurl}           % Not needed if you use pdflatex only.
\fi

\usepackage{multirow}
\let\cline\clineorig

\usepackage{amsmath,amssymb,amsfonts,bm}
\usepackage{algorithmic}
\usepackage{graphicx}
\usepackage{textcomp}
\usepackage{xcolor}
\usepackage{amsthm}

\newtheorem{hypothesis}{Hypothesis}

\usepackage{listings}
\usepackage{lstlangcoq}
\usepackage{xcolor}
\usepackage[T1]{fontenc}
\usepackage[scaled]{beramono}

\makeatletter
\AtBeginDocument{%
  \let\c@figure\c@lstlisting
  
  \let\ftype@lstlisting\ftype@figure % give the floats the same precedence
}
\makeatother

\newcommand\blfootnote[1]{%
  \begingroup
  \renewcommand\thefootnote{}\footnote{#1}%
  \addtocounter{footnote}{-1}%
  \endgroup
}

\title{Formal Verification of COO to CSR \\Sparse Matrix Conversion
  (Invited Paper)\blfootnote{This paper accompanies my keynote lecture ``Foundational End-to-end Verification of Numerical Programs'' at VSS 2025, the International Workshop on Verification of Scientific Software; and covers one of the results described
in that talk.}
}

\newcommand\comment[1]{}

\author{Andrew W. Appel
\institute{Princeton University}}
\def\titlerunning{Formal verification of COO to CSR Sparse Matrix Conversion}
\def\authorrunning{Andrew W. Appel}
\begin{document}

\maketitle
\thispagestyle{plain}\pagestyle{plain} % turn on page numbers for submission

\begin{abstract}
We describe a machine-checked correctness proof
  of a C program that converts a coordinate-form (COO) sparse matrix
  to a compressed-sparse-row (CSR) matrix.  The classic algorithm 
  (sort the COO entries in lexicographic order by row,column;
  fill in the CSR arrays left to right) is concise
  but has rather intricate invariants.  We illustrate a bottom-up
  methodology for deriving the invariants from the program.
\end{abstract}

%\begin{IEEEkeywords}
%numerical methods, formal verification, sparse matrix
%program verification
%\end{IEEEkeywords}

\section{Introduction}
We will describe the machine-checked correctness proof of a
C program that converts a Coordinate-form sparse matrix
(COO) into a Compressed Sparse Row matrix (CSR):
\begin{lstlisting}[language=C]
struct csr_matrix *coo_to_csr_matrix(struct coo_matrix *p) { . . . }
\end{lstlisting}
The program itself is given in \autoref{lst:coo-csr-conv}; it implements
an algorithm (presumably) known for many decades\footnote{Barret \emph{et al.}
\cite{barret94:templates}
describe CSR (which they called ``compressed row storage'')
but do not mention COO nor the method of constructing CSR matrices---however,
if CSR matrices were in use then there must have been a method of constructing them.}
This paper is meant as
a tutorial on the methodology for approaching the specification and
proof of numerical algorithms involving both data structures and
approximations (keeping the data-structure reasoning
separate from the approximation reasoning), and a demonstration
of a technique for deriving loop invariants from the
properties they must satisfy in a Hoare logic proof.

Sparse matrix-vector multiplication is a fundamental operation in
numerical methods, and takes time proportional to the number
of nonzero entries in the matrix---which can be much smaller than
the size of the corresponding dense matrix.
Depending on the structure of sparsity in the matrix,
and on whether the multiplication is of the form $Ax$ or $x^{T}A$,
different sparse representations may be appropriate.  For the very
common case of unstructured sparsity and multiplication $Ax$, the
Compressed Sparse Row (CSR) representation is useful
\cite[\S4.3.1]{barret94:templates}.

\begin{figure}
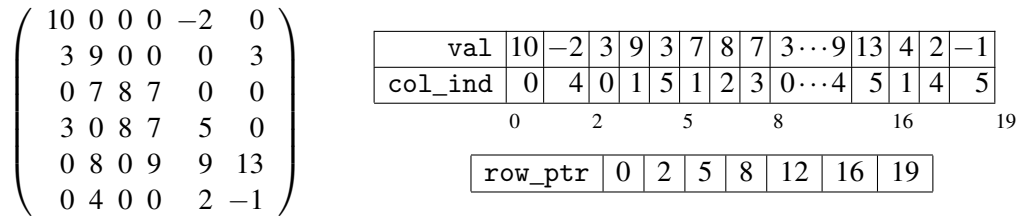

\(
\left(
\begin{array}{r@{~~}r@{~~}r@{~~}r@{~~}r@{~~}r}
  10 & 0 & 0 & 0 & -2 & 0 \\
  3 & 9 & 0 & 0 & 0 & 3 \\
  0 & 7 & 8 & 7 & 0 & 0 \\
  3 & 0 & 8 & 7 & 5 & 0 \\
  0 & 8 & 0 & 9 & 9 & 13 \\
  0 & 4 & 0 & 0 & 2 & -1 \\
\end{array}
\right)
\)\qquad\(
%\genfrac{}{}{0pt}{2} if using amsmath, otherwise:
\begin{array}{c}
\begin{array}{|r|@{\hspace{1pt}}r@{\hspace{2pt}}|@{\hspace{1pt}}r@{\hspace{2pt}}|@{\hspace{1pt}}r@{\hspace{2pt}}|@{\hspace{1pt}}r@{\hspace{2pt}}|@{\hspace{1pt}}r@{\hspace{2pt}}|@{\hspace{1pt}}r@{\hspace{2pt}}|@{\hspace{1pt}}r@{\hspace{2pt}}|@{\hspace{1pt}}r@{\hspace{2pt}}|@{\hspace{1pt}}r@{\hspace{2pt}}|@{\hspace{1pt}}r@{\hspace{2pt}}|@{\hspace{1pt}}r@{\hspace{2pt}}|@{\hspace{1pt}}r@{\hspace{2pt}}|@{\hspace{1pt}}r@{\hspace{2pt}}|@{}l@{}}
  \cline{1-14}
\mathtt{val} & 10 & -2 & ~3 & ~9 & ~3 & ~7 & ~8 & ~7 & ~3 \cdots9 & 13 & ~4 & ~2 &
-1 \\
  \cline{1-14}
\mathtt{col\_ind} & 0 & 4 & 0 & 1 & 5 & 1 & 2 & 3 & 0 \cdots 4 &5 & 1 & 4 & 5 \\
  \cline{1-14}
  \multicolumn{1}{l}{} &
  \multicolumn{2}{@{}l}{\mbox{\scriptsize 0}} &
  \multicolumn{3}{@{}l}{\mbox{\scriptsize{2}}} &
  \multicolumn{3}{@{}l}{\mbox{\scriptsize{5}}} &
  \multicolumn{2}{@{}l}{\mbox{\scriptsize{8}}} &
  \multicolumn{3}{@{}l}{\mbox{\scriptsize{16}}} &
  \multicolumn{1}{@{}l}{\mbox{\scriptsize{19}}} \\[6pt]
\end{array}\\
\begin{array}{|r|r|r|r|r|r|r|r|}
  \hline
  \mathtt{row\_ptr} & 0 & 2 & 5 & 8 & 12 & 16 & 19 \\
  \hline
\end{array}
\end{array}
\)
\caption{An example of the three arrays (\lstinline{val},
\lstinline{col\_ind}, \lstinline{row\_ptr}) used to store a matrix
compressed sparse row (CSR) format.
From~Barret \emph{et al.} \cite{barret94:templates}, adjusted for 0-based array indexing.
}
\label{fig:CSR}
\end{figure}

The CSR
format stores the elements of a sparse $m\times n$ matrix $A$
using three one-dimensional arrays: a floating-point array
\lstinline{val} that stores the nonzero elements of $A$, an
integer array \lstinline{col$\_$ind} that stores the column
indices of the elements, and an integer
array \lstinline{row$\_$ptr} that stores the locations in the array
\lstinline{col$\_$ind} that start a row in $A$. 
\autoref{fig:CSR} shows an example.

\begin{figure}
%\genfrac{}{}{0pt}{2} if using amsmath, otherwise:
\[
\begin{array}{|r|@{\hspace{1pt}}r@{\hspace{2pt}}|@{\hspace{1pt}}r@{\hspace{2pt}}|@{\hspace{1pt}}r@{\hspace{2pt}}|@{\hspace{1pt}}r@{\hspace{2pt}}|@{\hspace{1pt}}r@{\hspace{2pt}}|@{\hspace{1pt}}r@{\hspace{2pt}}|@{\hspace{1pt}}r@{\hspace{2pt}}|@{\hspace{1pt}}r@{\hspace{2pt}}|@{\hspace{1pt}}r@{\hspace{2pt}}|@{\hspace{1pt}}r@{\hspace{2pt}}|@{\hspace{1pt}}r@{\hspace{2pt}}|@{\hspace{1pt}}r@{\hspace{2pt}}|@{\hspace{1pt}}r@{\hspace{2pt}}|@{\hspace{1pt}}r@{\hspace{2pt}}|@{}l@{}}
  \cline{1-15}
\mathtt{row\_ind} & 0 & 0 &  0 & 1 & 1 & 1 & 2 & 2 & 2 & 2 \cdots 4 &  4 & 5 & 5 & 5 \\
\mathtt{col\_ind} & 0 & 0 &  4 & 0 & 1 & 5 & 1 & 2 & 3 & 5 \cdots 4 &  5 & 1 & 4 & 5 \\
\mathtt{val}      & 7 & 3 & -2 & 3 & 9 & 3 & 7 & 8 & 7 & 3 \cdots 9 & 13 & 4 & 2 & -1 \\
  \cline{1-15}
\end{array}
\]
\caption{A coordinate-form representation of the
matrix from \autoref{fig:CSR}}
\label{fig:COO}
\end{figure}

\begin{figure}
\begin{minipage}{3in}
\begin{lstlisting}[language=C,captionpos=b,caption= CSR struct in C.,label=lst:csrstruct]
struct csr_matrix {
     double *val; 
     unsigned *col_ind, *row_ptr;
     unsigned rows, cols; 
};
\end{lstlisting}
\end{minipage}\quad
\begin{minipage}{3in}
\begin{lstlisting}[language=C,captionpos=b,caption= COO struct in C.,label=lst:coostruct]
struct coo_matrix {
  unsigned *row_ind, *col_ind;
  double *val;
  unsigned n, maxn, rows, cols;
};
\end{lstlisting}
\end{minipage}
\end{figure}

But before a CSR matrix is used for multiplications, it must be \emph{built}.
One does not first build the dense matrix and from it extract the sparse
matrix, as that would be quite inefficient.  In a typical application
scientific/engineering application that gives rise to a sparse
matrix, one first translates the problem into a
a set of triples (row,column,value), that is, a \emph{coordinate form
sparse matrix} (COO matrix).  

A COO matrix has dimensions \lstinline{rows$\times$cols};
there are $n$ coordinate triples
\lstinline{row_ind[$k$]}, \lstinline{col_ind[$k$]}, \lstinline{val[$k$]}
for $0\le k < n$.
Each of those arrays has size \lstinline{maxn $\ge n$} to allow
for additional entries in the future.

A COO matrix may have more than one entry at the same (row,col).
If entry $k$ is $(i,j,x)$ and entry $k'$ is $(i,j,y)$,
the matrix this represents has value $x+y$ at position $(i,j)$---or $x+y$ plus other entries $(i,j,\_)$.  We call those \emph{duplicates.}
For example, the COO in \autoref{fig:COO} represents the same matrix
as the CSR in \autoref{fig:CSR}; there is a duplicate $7+3$ at position $(0,0)$.
Although the entries in \autoref{fig:COO} are sorted, generally
the entries in a COO matrix can be arranged in any order.

\begin{figure}
\begin{minipage}[b]{3in}
\caption{Mesh arising from an irregular finite-element problem.
{\footnotesize From \url{https://en.wikipedia.org/wiki/Mesh_generation}}
}
\label{fig:mesh}
\end{minipage}\quad
\begin{minipage}[b]{3.5in}
\includegraphics[scale=0.7]{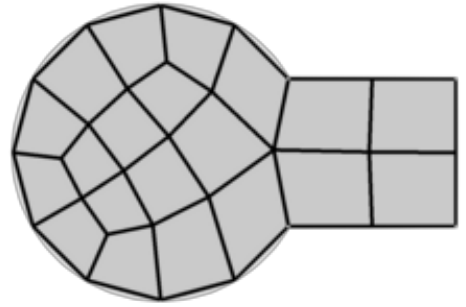}
\vspace*{-10pt}
\end{minipage}
\end{figure}

Duplicate entries arise naturally in scientific problems.
For example, in a finite-element analysis of a mesh such
as \autoref{fig:mesh},
each interior vertex adjacent to $d$ regions (elements) contributes 
$d$ values to the list of coordinate tuples.

COO matrices are not very efficient for performing matrix multiplication
(though perhaps better than dense matrices).  Their primary purpose is
as an intermediate representation for building another
form of sparse representation, such as CSR.

The algorithm for converting a COO to a CSR matrix is well known:  First, sort
the tuples by the lexicographic order of row then column.
Then process the sorted tuples in order, adding duplicates together
as they are seen; each of the three arrays of the CSR representation
can be filled in left-to-right.

\autoref{lst:coo-csr-conv} shows the program.  I have never run this
program, I have only proved it correct.  But unlike
Knuth,\footnote{``Beware of bugs in the above code; I have only proved
it correct, not tried it.''  Donald E. Knuth, Notes on the van Emde
Boas construction of priority queues, March 1977. https://staff.fnwi.uva.nl/p.vanemdeboas/knuthnote.pdf}
I can be confident
that it works, because (almost 50 years later) we have foundationally
sound machine-checked program logics.

The Coq proofs corresponding to this draft of the paper may be found at
{\footnotesize\url{https://github.com/VeriNum/iterative_methods/tree/8de4d78c78f92f280581253af49c2309fb95b2bc/sparse}}

\section{The specification}
In many application domains, we can prove that a program computes
exactly the right answer.  However, numerical algorithms generally
compute \emph{approximations}, with inaccuracies arising from
discretization and from floating-point round-off.  Thus the specification
of a numerical program will typically bound the distance between the
computed answer and the true solution to some mathematical equation.
In COO-to-CSR the approximation comes from floating-point
round-off when summing duplicate
elements.

Reasoning about data structures in C programs (or other low-level
imperative programming languages) is difficult enough
without having to simultaneously
reason about the composition of approximations.  So we separate
the reasoning into: first, the C program \emph{exactly}
implements
a low-level specification; and second, that low-level spec approximates the
solution to a mathematical equation within specified bounds.

We will use the Verifiable Software Toolchain to specify and prove correctness of
the C program.  In VST program logic, called \emph{Verifiable C},
we distinguish between mathematical objects and the data structures that represent them.
Verifiable C is embedded in the Coq logic, and
we can describe mathematical objects directly using Coq types and values.  For example, we consider ``COO matrix'' be a mathematical value whose
Coq type is,

\begin{lstlisting}
Record coo_matrix (t: type) := {
  coo_rows: $\mathbb{Z}$;
  coo_cols: $\mathbb{Z}$;
  coo_entries: list ($\mathbb{Z}$ * $\mathbb{Z}$ * ftype t)
}.
\end{lstlisting}
Here, \lstinline{type} means ``floating-point format,'' such as IEEE double precision;
for generality, our matrix types and our algorithms are parameterized over formats.
If \lstinline{t} is a \lstinline{type}, then
\lstinline{ftype t} is the Coq \lstinline{Type} of values in that floating-point format \cite{appel24:vcfloat}.

Therefore, a \lstinline{coo_matrix} is a tuple $(R,C,E)$ where
the matrix is supposed to have dimension $R\times C$
and $E$ is a list of 3-tuples.

A program might represent COO matrices with an array of 3-element records,
or Fortran-like with three separate arrays.  Either way, we can describe the relation
between the mathematical value of type \lstinline{coo_matrix(Tdouble)}
and our data-structure layout.  In Verifiable C such a relation is a
\emph{memory predicate} or \lstinline{mpred}.  So, for example, in our specification we have,

\begin{lstlisting}
Definition coo_rep (sh: share) (coo: coo_matrix Tdouble) (p: val) : mpred :=  . . .
\end{lstlisting}
The permission-share \lstinline{sh} tells whether the data structure has
read permission or write permission (etc.) and henceforth we will ignore or omit
permission shares.  (They can be quite useful when describing shared-memory
parallel programs, but that is not the focus here.)
We can read \lstinline{coo_rep} as saying, the COO matrix \lstinline{coo}
is laid out in memory as a pointer data structure rooted at address \lstinline{p}.

\pagebreak
Similarly, we can describe the mathematical type of CSR matrices:

\begin{lstlisting}
Record csr_matrix {t: type} := {
  csr_cols: $\mathbb{Z}$;
  csr_vals: list (ftype t);
  csr_col_ind: list $\mathbb{Z}$;
  csr_row_ptr: list $\mathbb{Z}$;
  csr_rows: $\mathbb{Z}$ := Zlength (csr_row_ptr) - 1
}.
\end{lstlisting}
That is, a CSR matrix is a 4-tuple (cols,vals,col\_ind,row\_ptr)
where cols is the number of columns in the matrix,
vals is a sequence of floating-point values in format $t$,
col\_ind is a sequence of integers,
and row\_ptr is another sequence of integers;
and where the number of rows is one less than the length of the row\_ptr sequence.
See \autoref{fig:CSR} for an example.

One can imagine various data structures with which a C program could represent this structure
in memory, but having chosen such a data structure, we can specify it with an \lstinline{mpred}
relation:
\begin{lstlisting}
Definition csr_rep (sh: share) (csr: csr_matrix Tdouble) (p: val) : mpred :=  . . .
\end{lstlisting}

A COO matrix represents a mathematical matrix, and we can state this with a mathematical relation:

\begin{lstlisting}
Definition coo_to_matrix {t} (coo: coo_matrix t) (m: matrix t) : Prop := . . .

Definition csr_to_matrix {t} (csr: csr_matrix t) (m: matrix t) : Prop := . . .
\end{lstlisting}

So, for example, the relation \lstinline{csr_to_matrix} holds between
the two mathematical objects depicted in \autoref{fig:CSR};
and the relation \lstinline{coo_to_matrix} holds between
the COO shown in \autoref{fig:COO} and the matrix in \autoref{fig:CSR}.

Having defined all these relations, we can now specify what it would mean for
the C function \lstinline{coo_to_csr_matrix} to be correct.
Our program takes as input the address $p$ where a COO matrix is stored, and
returns the address at which a corresponding CSR matrix is stored.

One might think the specification of this program is as follows:

\begin{itemize}
\item  Let $\mathit{coo}$ be a COO sparse matrix representation,
\item that is laid out in memory at address $p$ ~ (i.e., \lstinline{coo_rep $\mathit{coo}$ $p$});
\item let $M$ be the matrix that $\bar{M}$  represents (i.e., \lstinline{coo_to_matrix $\mathit{coo}$ $M$});
\item then there exists a CSR sparse matrix $\mathit{csr}$,
\item that represents matrix $M$,
\item and that is laid out in memory at address $q$ ~ (i.e., \lstinline{csr_rep $\mathit{csr}$ $q$});
\item and $q$ is returned.
\end{itemize}

This is almost right.  The problem is the use of the definite article,
``let $M$ be \emph{the} matrix.''  A floating-point COO matrix does not represent a unique
mathematical matrix, because of the floating-point addition needed when duplicate
entries are combined.  The duplicate entries may be added together in any order;
and while addition in the reals is associative, addition in the floats is not.
The possible values of the resulting matrix-entry will all be similar to each other,
but not exactly the same.
That is, \lstinline{coo_to_matrix} is a relation but it is not a function.
We can say, ``let $M$ be \emph{a} matrix that $\bar{M}$ represents.''

To define the relation between $\bar{M}$ and $M$, we first need a relation defining the floating point sum, \emph{in any order and any tree of associativity,} of a set of floating-point values:

\begin{lstlisting}
Inductive sum_any {t}: list (ftype t) -> ftype t -> Prop :=
| Sum_Any_0: sum_any nil (Zconst t 0)
| Sum_Any_1: forall x, sum_any [x] x
| Sum_Any_split: forall al bl a b, sum_any al a -> sum_any bl b ->
                  sum_any (al++bl) (BPLUS a b)
| Sum_Any_perm: forall al bl s, Permutation al bl -> sum_any al s -> sum_any bl s.
\end{lstlisting}
\noindent That is, let $t$ be a floating-point format, so the Coq type
\lstinline{ftype t} contains floating-point values in that format.
The relation \lstinline{sum_any $v$ $s$}, where the $t$ argument is implicit,
says that the list-of-floats $v$  relates to float $s$ as follows:
\begin{itemize}
\item If $v$ is the empty list then $s=0$.
\item If $v$ is the singleton list containing $x$ then $s=x$.
\item If $v$ is the concatentation of lists $\mathit{al}$ and $\mathit{bl}$,
such that $\mathit{al}$ sums (in any order) to $a$ and $\mathit{bl}$ to $b$,
then $v$ relates to the floating-point sum of $a$ and $b$,
written as \lstinline{BPLUS $a$ $b$}.
\item If $\mathit{al}$ is a permutation of $\mathit{bl}$,
and $\mathit{al}$ relates to $s$, then $\mathit{bl}$ also relates to $s$.
\end{itemize}

By this definition, \lstinline{sum_any} relates $\mathit{al}$ to any
floating-point value that relates from combining the elements of $s$,
once each, in an arbitrary tree of additions.  There are theorems in numerical analysis that bound the distance between any such $s$ and the real number
that one would obtain by adding in arbitrary-precision arithmetic without rounding.  Such theorems will be useful later, but we don't need them now;
we can first specify correctness of COO-to-CSR conversion using
\lstinline{sum_any}, then show that this specification implies
the desired accuracy bound.

We can therefore write the definition of \lstinline{coo_to_matrix}---that a COO matrix \lstinline{coo} represents a mathematical matrix \lstinline{m}:
\begin{lstlisting}[language=Coq,captionpos=b,caption=coo\_to\_matrix,label=lst:coo-to-matrix]
Definition coo_to_matrix {t: type} (coo: coo_matrix t) (m: matrix t) : Prop :=
  coo_rows coo = matrix_rows m /\
  matrix_cols m (coo_cols coo) /\
   forall i, 0 $\le$ i < coo_rows coo ->
    forall j, 0 $\le$ j < coo_cols coo -> 
     sum_any (map snd (filter (coord_eqb (i,j) oo fst) (coo_entries coo)))
          (matrix_index m (Z.to_nat i) (Z.to_nat j)).
\end{lstlisting}
At index $(i,j)$ we filter all the entries of the COO matrix whose
coordinates are equal to $(i,j)$,
take their values (with \lstinline{map snd}),
and use  \lstinline{sum_any} to relate that to $m_{ij}$.

\begin{lstlisting}[language=Coq,float,captionpos=b,caption={Specification (simplified) of the C function.  The full specification describes how this function has access to the memory allocator (malloc/free), needed to allocate space for the new CSR matrix.},label=lst:funspec]
Definition coo_to_csr_matrix_spec [simplified] :=
 DECLARE _coo_to_csr_matrix
 WITH coo: coo_matrix Tdouble,
 PRE [ tptr (Tstruct _coo_matrix noattr) ]
    PROP(coo_matrix_wellformed coo)
    PARAMS( p )
    SEP (coo_rep coo p)
 POST [ tptr (Tstruct _csr_matrix noattr) ]
   EX coo': coo_matrix Tdouble, EX csr: csr_matrix Tdouble, 
       EX m: matrix Tdouble, EX q: val,
    PROP(coo_matrix_equiv coo coo';
         coo_to_matrix coo m; csr_to_matrix csr m)
    RETURN( q )
    SEP (coo_rep coo' p; csr_rep csr q).
\end{lstlisting}
\autoref{lst:funspec} gives the ``funspec'' (VST function specification)
for the \lstinline{coo_to_csr_matrix} function.
We will summarize here what it means; for more explanation,
see \cite{vst-floyd} or \cite{VerifiableCmanual}.
The lines of the funspec are as follows:
\begin{description}
\item[DECLARE:] The C name of the function is coo\_to\_csr\_matrix.
\item[WITH] quantifies over the mathematical variables
that will be used in the precondition (and perhaps postcondition):
the COO matrix \lstinline{coo} and address $p$.
\item[PRE] begins the precondition of the function, which takes a C parameter
whose C type is ``pointer to struct csr\_matrix.''
\item[PROP:]  The precondition has three parts, this first of which
contains mathematical propositions that must hold of the WITH-bound
variables; in this case, that the COO matrix is well-formed (for example,
that each entry's row- and column-index are within the bounds of the
matrix dimensions; see \autoref{lst:wellformed}).
\item[PARAM:] The value of the C-language parameter is
\lstinline{p}.  We distinguish between mathematical values bound in the
logic (such as \lstinline{p}) from C-language identifiers that
stand for C variables (such as \lstinline{_p}).  In this case,
the C variable \lstinline{_p} \emph{contains} the value
\lstinline{p} upon entry to the function.
\item[SEP:] This clause contains spatial conjuncts in separation logic,
that is, it describes data structures in memory.  In this case,
that there is a representation of the \lstinline{coo} matrix
at address $p$.
\item[POST:] The postcondition starts by giving the C type of the
return value of the function; in this case, pointer to 
\lstinline{struct csr_matrix}.
\item[EX:] This postcondition existentially quantifies four mathematical
quantities: a COO matrix \lstinline{coo'}, a CSR matrix \lstinline{csr},
a mathematical matrix \lstinline{m}, and an address \lstinline{q}.
\item[PROP:]
Three propositions are asserted about how the variables are related:
\lstinline{coo'} is equivalent to \lstinline{coo},
\lstinline{coo} represents the matrix \lstinline{m},
and \lstinline{csr} also represents \lstinline{m}.
The postcondition doesn't say so explicitly, but
in a typical implementation
the difference between \lstinline{coo} and \lstinline{coo'} is that
the entries (coordinate-tuples) of \lstinline{coo'} are sorted in order.
\item[RETURN:] address \lstinline{q} is the value returned by this function.
\item[SEP:]  In memory when the function returns are the modified COO
matrix \lstinline{coo'} (at the same address \lstinline{p}) and a
CSR matrix at newly allocated memory address \lstinline{q}.
\end{description}

\begin{lstlisting}[language=Coq,float,
captionpos=b,caption= Wellformedness of COO matrices,label=lst:wellformed]
Definition coo_matrix_wellformed {t} (coo: coo_matrix t) :=
 (0 $\le$ coo_rows coo /\ 0 $\le$ coo_cols coo)
 /\ Forall (fun e => 0 $\le$ fst (fst e) < coo_rows coo /\ 0 $\le$ snd (fst e) < coo_cols coo)
      (coo_entries coo).
\end{lstlisting}

\begin{lstlisting}[language=C,float,numbers=left,numberstyle=\tiny,
captionpos=b,caption= COO-to-CSR conversion,label=lst:coo-csr-conv]
unsigned coo_count (struct coo_matrix *p) {
  unsigned i, n = p$~$->n;
  if (n==0) return 0;
  unsigned count=1;
  for (i=1; i<n; i++)
    if (p$~$->row_ind[i-1]!=p$~$->row_ind[i] || p$~$->col_ind[i-1]!=p$~$->col_ind[i])
      count++;
  return count;
}

struct csr_matrix *coo_to_csr_matrix(struct coo_matrix *p) {
  struct csr_matrix *q;
  unsigned count, i, r,c, ri, ci, cols, k, l, rows;
  unsigned *col_ind, *row_ptr, *prow_ind, *pcol_ind;
  double x, *val, *pval;
  unsigned n = p$~$->n;
  coo_quicksort(p, 0, n);
  k = coo_count(p);
  rows = p$~$->rows; prow_ind=p$~$->row_ind; pcol_ind=p$~$->col_ind;
  pval = p$~$->val;
  q = surely_malloc(sizeof(struct csr_matrix));
  val = surely_malloc(k * sizeof(double));
  col_ind = surely_malloc(k * sizeof(unsigned));
  row_ptr = surely_malloc ((rows+1) * sizeof(unsigned));
  r=-1;
  c=0; /* this line is unnecessary but simplifies the proof */
  l=0;                 /* partial_csr_0 */
  for (i=0; i<n; i++) {
    ri = prow_ind[i];  ci = pcol_ind[i]; x = pval[i];
    if (ri==r)
      if (ci==c)
	val[l-1] += x; /* partial_CSR_duplicate */
      else { c=ci; col_ind[l]=ci; val[l]=x; l++;}  /* partial_CSR_newcol */
    else {
      while (r+1<=ri) row_ptr[++r]=l;              /* partial_CSR_skiprow */
      c=ci; col_ind[l]=ci; val[l]=x; l++;          /* partial_CSR_newrow */
  } }
  cols = p$~$->cols;
  while (r+1<=rows) row_ptr[++r]=l;  /* partial_CSR_lastrows */
  q$~$->val = val;  q$~$->col_ind = col_ind; q$~$->row_ptr = row_ptr;
  q$~$->rows = rows;  q$~$->cols = cols;
  return q;          /* partial_CSR_properties */
}
\end{lstlisting}

\subsection{The minimum you need to inspect}
When one inspects a machine-checked proof, it is not so important
to check the proof itself---the machine has done that---but
it is critical to check the theorem-statement.  For if
it's the wrong theorem, it will not serve the intended purpose.
Furthermore, every definition that's referenced (even indirectly)
from the theorem-statement is part of this ``trusted base.''
So let us review what goes into the theorem we have stated.
\begin{itemize}
\item \lstinline{coo_to_csr_matrix_spec} (\autoref{lst:funspec}), 14 lines.
\item \lstinline{coo_to_matrix} (\autoref{lst:coo-to-matrix}), 7 lines.
\item \lstinline{coo_rep} (not shown), 15 lines.
\item \lstinline{coo_matrix_equiv} (not shown), 3 lines.
\end{itemize}
Although \lstinline{csr_to_matrix} and \lstinline{csr_rep}
are mentioned in \lstinline{coo_to_csr_matrix_spec},
it is not strictly necessary to inspect them.  That's because
we have elsewhere \cite{kellison23:laproof}
proved the theorem that our C-language sparse matrix-vector multiply
function correctly multiplies a CSR matrix by a vector to
produce the correct result that the mathematical matrix $M$ would
produce (based on the same definitions of \lstinline{csr_to_matrix}
and \lstinline{csr_rep}).\footnote{The definition of \lstinline{csr_rep}
 in the LAProof paper \cite{kellison23:laproof}
encompasses what, it this paper, we have separated
 into \lstinline{csr_to_matrix} and \lstinline{csr_rep}.}
Therefore, you can treat these definitions as an abstract data type;
no matter what they are, you will get the intended result when multiplying
by the CSR matrix that our function produces.

In the remainder of this paper I will present dozens of lines of
definitions and lemma-statements.  All of those should be treated
as part of the \emph{proof} of the main theorem, stated above;
and this proof has been checked by the Coq kernel.

\subsection{An alternate specification}
Some computational scientists, since they cannot normally achieve \emph{formal verification
of correctness and accuracy} such as we do here, rely upon
other means of validation and verification.  One of those is \emph{regression testing;}
that is, testing whether a change to the program has introduced a new bug
by comparing the output of the changed program to the output of the original program,
or examining the output of specific test cases.
For such testing, it is helpful if the program has
the property of \emph{bit-for-bit reproducibility}.  The COO-to-CSR program
as I have specified it above does not have that property, because it
permits any order of summation of duplicate elements, and in floating point
those values may be very slightly different.

One can write a stronger specification, that gives bit-for-bit reproducibility,
by insisting that the duplicate elements be added in left-to-right order
of their appearence in the original (unsorted) COO entry list.  Then one can
easily show that the stronger specification entails the specification I've given above.
The program I have verified in this paper does not have that property.

A separate issue is that of \emph{reassembly}.  Suppose one converts a COO matrix to CSR,
and then one has a new (i.e., modified) COO matrix with the same sparsity structure
but different values.  Then the structure of the CSR matrix will be unchanged,
and it should be possible to update just the \lstinline{val} array of the CSR matrix,
slightly more efficiently than building a CSR matrix from scratch.
To specify this program, and at the same time achieve bit-for-bit reproducibility,
we would divide the \lstinline{coo_to_csr_matrix} into two parts:
\begin{enumerate}
\item Build the structure of the CSR matrix (\lstinline{row_ptr} and \lstinline{col_ind})
\item Fill in the values.
\end{enumerate}
Such a program would be similar in many ways to the program I verify here,
and would satisfy very similar loop invariants.  However, we leave this for future work.

\section{The low-level proof}
It is useful to stratify a proof to separate the low-level details
of C programs and data structures from high-level concerns about
algorithms.  Often this is done by stating an algorithm as a pure
functional program in Coq, proving that the C program refines the
algorithm, then proving that the algorithm is correct.  Here we will
take a different approach: we will prove that the C program
is correct provided a certain relation exists with certain properties;
then we will give a model of the relation.  

Let us examine the program (\autoref{lst:coo-csr-conv}).
First (at line 17) the COO entries
are sorted in place. Then the \lstinline{coo_count} function
scans the entries in order, counting how many
\emph{distinct} (row,column) coordinate-pairs are in
the list.  This allows (at lines 22--23) the allocation of
CSR arrays of exactly the right size.\footnote{\lstinline{surely_malloc}
calls \lstinline{malloc} but then, if malloc returns NULL indicating
it failed to allocate memory, \lstinline{surely_malloc} aborts
the program.}

The main loop begins at line 28, traversing all the entries of the
COO matrix in order.  In examining each entry $(r_i,c_i,x_i)$
we need to check whether it is a duplicate of some previous entry
$(r,c,y)$, for which purpose the local variables record
the previous entry's $r$ and $c$. The previous $y$ need not be remembered because it has already been
added into the appropriate spot of the \lstinline{val} array.

We use unsigned integers for row and column numbers, in part because
that allows a greater range of indexes, i.e., bigger matrices.

To ensure that the first entry is not treated as a duplicate,
(at line 25) we initialize $r$ to one less than the smallest possible row number; that is, $-1$.
Since we are using unsigned (i.e., modulo $2^w$) arithmetic, this is really
$2^w-1$, so we must be careful.  The program \emph{is} careful, but the version before I 
debugged it (by proving it correct) was not as careful.  For example,
\lstinline{while (r+1<=ri)} is correct at line 35,
but \lstinline{while (r<ri)} would be wrong.  One can learn this
(as I did) from the fact that VST's proof system, since it is sound,
won't permit a proof of an incorrect program.

In the traversal, the variable \lstinline{l} tracks the number
of \emph{distinct} coordinate-pairs seen so far; this indicates
the spot in \lstinline{val} array (and in the \lstinline{col_ind} array)
to be filled in next.

\begin{figure}[h]
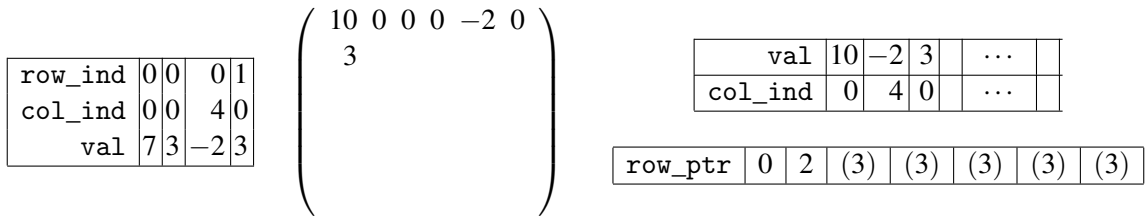

\(
\begin{array}{|r|@{\hspace{1pt}}r@{\hspace{2pt}}|@{\hspace{1pt}}r@{\hspace{2pt}}|@{\hspace{1pt}}r@{\hspace{2pt}}|@{\hspace{1pt}}r@{\hspace{2pt}}|@{}l@{}}
  \cline{1-5}
\mathtt{row\_ind} & 0 & 0 &  0 & 1 \\
\mathtt{col\_ind} & 0 & 0 &  4 & 0 \\
\mathtt{val}      & 7 & 3 & -2 & 3 \\
  \cline{1-5}
\end{array}
\)
\quad
\(
\left(
\begin{array}{r@{~~}r@{~~}r@{~~}r@{~~}r@{~~}r}
  10 & 0 & 0 & 0 & -2 & 0 \\
  3 &  &  &  &  &  \\
   &  &  &  &  &  \\
   &  &  &  &  &  \\
   &  &  &  &  &  \\
   &  &  &  &  & \\
\end{array}
\right)
\)
\quad
\(
%\genfrac{}{}{0pt}{2} if using amsmath, otherwise:
\begin{array}{c}
\begin{array}{|r|@{\hspace{1pt}}r@{\hspace{2pt}}|@{\hspace{1pt}}r@{\hspace{2pt}}|@{\hspace{1pt}}r@{\hspace{2pt}}|@{\hspace{1pt}}r@{\hspace{2pt}}|@{\hspace{1pt}}r@{\hspace{2pt}}|@{\hspace{1pt}}r@{\hspace{2pt}}|@{\hspace{1pt}}r@{\hspace{2pt}}|@{}l@{}}
  \cline{1-7}
\mathtt{val} & 10 & -2 & ~3 & ~~  &  \cdots ~~  & ~~  \\
  \cline{1-7}
\mathtt{col\_ind} & 0 & 4 & 0 &  &  ~~ \cdots ~~ &  \\
  \cline{1-7}
\end{array}\\
~\\
\begin{array}{|r|r|r|r|r|r|r|r|}
  \hline
  \mathtt{row\_ptr} & 0 & 2 & (3) & (3) & (3) & (3) & (3) \\
  \hline
\end{array}
\end{array}
\)
\caption{A partially completed CSR representation corresponding to the first four
  entries in the COO matrix of \autoref{fig:COO}.  The (3) values
  in \lstinline{row_ptr} are not actually stored.}
\label{fig:CSR-partial}
\end{figure}

It is straightforward to do a Hoare-logic forward proof in VST
to reach line 28.  But then we will need a loop invariant.
In fact, there are 
three loops (at lines 28, 35, and 39); we will derive all of their
invariants together.
At any point during the execution of the loop(s), the program
has built a partial CSR matrix that represents, more or less,
the first $i$ entries of of the (sorted) COO matrix.
This (partial) CSR matrix is represented in the arrays
\lstinline{val}, \lstinline{col_ind}, and \lstinline{row_ptr}.
\autoref{fig:CSR-partial} illustrates such a configuration.

Before we try to find a model (a definition) for this relation,
let us examine the properties such a relation must have.  At certain
points in the program, the values of variables $i$ and $r$ or the contents of \lstinline{rowptr, colind, val}
are changed.  These points are labeled in \autoref{lst:coo-csr-conv}
with comments of the form \lstinline{/* partial_CSR_... */}.

At each such point, we can examine the VST assertion that characterizes
the current program state, along with the assertion required
by the next iteration of the loop, to derive an axiom that
the \lstinline{partial_CSR} relation must satisfy.

\newsavebox\foo
\sbox{\foo}{\lstinline{partial_CSR i r coo rowptr colind val}}
\begin{hypothesis}
\label{hypo}
There is a relation ~\usebox\foo{}~
such that the arrays \lstinline{rowptr, colind, val}
represent the COO matrix entries up to the $i$th entry
and the $r$th row of the matrix;
and furthermore \lstinline{partial_CSR} satisfies all the lemma-statements
of \autoref{lst:axioms} and \autoref{lst:axioms2}.
\end{hypothesis}

\begin{lstlisting}[language=Coq,float,
captionpos=b,caption=Axioms for \lstinline{partial_CSR},label=lst:axioms]
Definition cd {t} (coo: coo_matrix t) := 
    count_distinct (coo_entries coo).
Definition cd_upto {t} i (coo: coo_matrix t) := 
    count_distinct (sublist 0 i (coo_entries coo)).

Lemma partial_CSR_0:
  forall (coo: coo_matrix Tdouble), 
  coo_matrix_wellformed coo -> 
  sorted coord_le (coo_entries coo) -> 
  cd coo $\le$ Int.max_unsigned -> 
  partial_CSR 0 (-1) coo 
      (Zrepeat Vundef (coo_rows coo + 1)) 
      (Zrepeat Vundef (cd coo)) (Zrepeat Vundef (cd coo)).

Lemma partial_CSR_duplicate:
    forall h r coo (f: ftype Tdouble) ROWPTR COLIND VAL,
    0 < h < Zlength (coo_entries coo) ->
    fst (Znth (h-1) (coo_entries coo)) = fst (Znth h (coo_entries coo)) ->
    r = fst (fst (Znth (h-1) (coo_entries coo))) ->
    Znth (cd_upto h coo - 1) VAL = Vfloat f ->
    partial_CSR h r coo ROWPTR COLIND VAL ->
    partial_CSR (h+1) r coo ROWPTR COLIND 
      (upd_Znth (cd_upto h coo - 1) VAL
            (Vfloat (Float.add f (snd (Znth h (coo_entries coo)))))).

Lemma partial_CSR_newcol:
   forall i r c x coo ROWPTR COLIND VAL,
   0 < i < Zlength (coo_entries coo) ->
   Znth i (coo_entries coo) = (r, c, x) ->
   r = fst (fst (Znth (i-1) (coo_entries coo))) ->
   c <> snd (fst (Znth (i-1) (coo_entries coo))) ->
   partial_CSR i r coo ROWPTR COLIND VAL ->
   partial_CSR (i + 1) r coo ROWPTR
     (upd_Znth (cd_upto i coo) COLIND (Vint (Int.repr c)))
     (upd_Znth (cd_upto i coo) VAL (Vfloat x)).

Lemma partial_CSR_skiprow:
   forall i r coo ROWPTR COLIND VAL,
   0 $\le$ i < Zlength (coo_entries coo) ->
   r $\le$ fst (fst (Znth i (coo_entries coo))) ->
   partial_CSR i (r-1) coo ROWPTR COLIND VAL ->
   partial_CSR i r coo (upd_Znth r ROWPTR (Vint (Int.repr (cd_upto i coo))))
     COLIND VAL.
\end{lstlisting}

\begin{lstlisting}[language=Coq,float,
captionpos=b,caption=More axioms for \lstinline{partial_CSR},label=lst:axioms2]
Lemma partial_CSR_newrow: 
    forall i r c x coo ROWPTR COLIND VAL,
    0 $\le$ i < Zlength (coo_entries coo) ->
    Znth i (coo_entries coo) = (r,c,x) ->
    (i <> 0 -> fst (fst (Znth (i - 1) (coo_entries coo))) <> r) ->
    partial_CSR i r coo ROWPTR COLIND VAL ->
    partial_CSR (i + 1) r coo ROWPTR
     (upd_Znth (cd_upto i coo) COLIND (Vint (Int.repr c)))
     (upd_Znth (cd_upto i coo) VAL (Vfloat x)).

Lemma partial_CSR_lastrows:
   forall r coo ROWPTR COLIND VAL,
    r $\le$ coo_rows coo ->
   partial_CSR (Zlength (coo_entries coo)) (r-1) coo ROWPTR COLIND VAL ->
   partial_CSR (Zlength (coo_entries coo)) r coo 
     (upd_Znth r ROWPTR (Vint (Int.repr (cd coo)))) COLIND VAL.

Lemma partial_CSR_properties:
  forall coo ROWPTR COLIND VAL,
    partial_CSR (Zlength (coo_entries coo)) (coo_rows coo) coo ROWPTR COLIND VAL ->
    exists (m: matrix Tdouble) (csr: csr_matrix Tdouble),
            csr_to_matrix csr m /\ coo_to_matrix coo m
            /\ coo_rows coo = matrix_rows m 
            /\ coo_cols coo = csr_cols csr 
            /\ map Vfloat (csr_vals csr) = VAL
            /\ Zlength (csr_col_ind csr) = cd coo
            /\ map Vint (map Int.repr (csr_col_ind csr)) = COLIND
            /\ map Vint (map Int.repr (csr_row_ptr csr)) = ROWPTR
            /\ Zlength (csr_vals csr) = cd coo.
\end{lstlisting}

Since the COO entries are sorted in order, first by row and then by
column, one might think that $i$ determines $r$ and that it is
unnecessary to include $r$ as a parameter of the relation.
But not quite.  The matrix might have some all-zero rows;
the purpose of the loops at lines
35 and 39 is to skip over those rows while recording information
about them into the \lstinline{row_ptr} array.  So the relation
must take note of both $i$ and $r$.

Even though we don't yet know the definition of this relation,
we can derive the axioms it must satisfy by attempting
a Hoare-logic proof of the program.  For the main loop,
we choose an invariant (\autoref{lst:inv}) that includes
\begin{lstlisting}
partial_CSR i r coo' ROWPTR COLIND VAL
\end{lstlisting}
where $i$ is the value of the loop iteration variable,
$r$ is the value local variable \lstinline{r},
\lstinline{coo'} is the COO matrix that results from
the quicksort called at line 17 (of \autoref{lst:coo-csr-conv}), and the remaining variables
are the current contents of the arrays
\lstinline{row_ptr}, \lstinline{col_ind}, and \lstinline{val}.
I have previously proved the correctness of a quicksort
implementation,\footnote{https://github.com/cverified/cbench-vst/tree/master/qsort, May 2019}
so I can use that function-specification
in the current proof to establish that the entries of
\lstinline{coo'} are sorted by the relation \lstinline{coord_le}.

Consider, then, the situation at line 28 just before
the loop starts.   We can
assume the COO matrix is well formed, that \lstinline{coo'}
is sorted, and the number of distinct elements in the
\lstinline{coo_entries} of \lstinline{coo'} is representable
as an unsigned integer, that is, \lstinline{$\le$ Int.max_unsigned}.
We have $r= -1$, $i=0$, and all three CSR arrays
are completely uninitialized.  We can summarize this with
a lemma called \lstinline{partial_CSR_0} (see \autoref{lst:axioms})
which concludes,
\begin{lstlisting}
 partial_CSR 0 (-1) coo 
      (Zrepeat Vundef (coo_rows coo + 1))    (* ROWPTR *)
      (Zrepeat Vundef (cd coo))              (* COLIND *)
      (Zrepeat Vundef (cd coo)).             (* VAL *)
\end{lstlisting}
That is, \lstinline{partial_CSR} relates the first 0 entries of
 the COO matrix
 to a partial CSR matrix filled up to rows=$-1$ such that the
 ROWPTR array is a sequence of (rows+1) undefined values,
 the COLIND array is a sequence of undefined values of length
 \lstinline{cd coo}, and the VAL array is
 a sequence of undefined (i.e., uninitialized) values of length
 \lstinline{cd coo}; where 
 \lstinline{cd coo} (``count distinct'') 
 is the number of distinct coordinate-pairs
 in the entry list of \lstinline{coo}.  We will need this many
 slots to be available for processing all the entries.

Consider another example: at line 33 we have determined that
the entry $(r_i,c_i,x)$ is not a duplicate, and we must create
a new column in the current row.  This corresponds to the lemma named
\lstinline{partial_CSR_newcol} (in \autoref{lst:axioms})
which concludes,
\begin{lstlisting}
 partial_CSR i r coo ROWPTR COLIND VAL ->
 partial_CSR (i + 1) r coo ROWPTR
     (upd_Znth (cd_upto i coo) COLIND (Vint (Int.repr c)))
     (upd_Znth (cd_upto i coo) VAL (Vfloat x)).
\end{lstlisting}
This says that if we increase $i$ to $i+1$,
and update the $k$th element of COLIND and ROW to $c$ and $x$
respectively (where $k$ is the number of distinct elements in
the first $i$ entries of \lstinline{coo}), then
the \lstinline{partial_CSR} entry still holds (provided the other
premises of \lstinline{partial_CSR_newcol}  are satisfied).
\lstinline{Int.repr} forces an mathematical integer into
a modulo-$2^w$ machine integer, and \lstinline{Vint} injects
that into the type of C values; \lstinline{Vfloat} injects
 an IEEE double-precision floating-point number into the C value type.

\autoref{lst:inv} shows the invariant for the main
loop.
At line 33 of \autoref{lst:coo-csr-conv} we must reestablish that invariant;
the lemma statement \lstinline{partial_CSR_newcol}
is derived by noticing what it will take to reestablish
this loop invariant.
The other lemma statements are derived analogously, as we
reestablish invariants at the places in the program annotated with
comments of the form \lstinline{/* partial_CSR_... */}.

Finally, once we reach the end of the loops (at line 42),
the relation \lstinline{partial_CSR} with $i=n$ and $\mathit{r=rows}$,
(which is to say, a completed CSR matrix)
should have properties strong enough to derive the desired
postcondition of the function; namely, those stated in
\lstinline{partial_CSR_properties}  (\autoref{lst:axioms2}).

Assuming the existence of a \lstinline{partial_CSR} relation
satisfying these axioms allows the C program to be proved correct
without too much fuss.
Now, if we can demonstrate that such a relation exists,
then the main theorem will be proved.

\begin{lstlisting}[language=Coq,float,
captionpos=b,caption={Main loop invariant, in VST notation.
The EX clauses bind existentially quantified variables.
The PROP clauses state propositions about those variables.
The LOCAL clauses give the values of C variables.
The SEP clauses give the separating conjuncts for
data structures in memory.  FRZL FR1 is essentially a ``frame''
of data that are untouched by (and irrelevant to) the loop.
},label=lst:inv]
EX i:$\mathbb{Z}$, EX l:$\mathbb{Z}$, EX r:$\mathbb{Z}$, EX c:$\mathbb{Z}$, 
EX ROWPTR: list val, EX COLIND: list val, EX VAL: list val,
  PROP(0$\le$l$\le$k; l$\le$i$\le$n; -1 $\le$ r < coo_rows coo'; 0 $\le$ c $\le$ coo_cols coo';
       partial_CSR i r coo' ROWPTR COLIND VAL;
       l = count_distinct (sublist 0 i (coo_entries coo'));
       l=0 -> r=-1;
       i<>0 -> r=(fst (fst (Znth (i-1) (coo_entries coo'))))%Z /\ c = snd (fst (Znth (i-1) (coo_entries coo')))) 
 LOCAL (temp _l (Vint (Int.repr l));
       temp _r (Vint (Int.repr r)); temp _c (Vint (Int.repr c));
       temp _row_ptr rowptr_p; temp _col_ind colind_p; temp _val val_p;
       temp _q q; temp _pval vp; temp _pcol_ind cp; temp _prow_ind rp;
       temp _rows (Vint (Int.repr (coo_rows coo')));
       temp _n (Vint (Int.repr n)); temp _p p)
  SEP(FRZL FR1;
      data_at Ews (tarray tuint (coo_rows coo' + 1)) ROWPTR rowptr_p;
      data_at Ews (tarray tuint k) COLIND colind_p; 
      data_at Ews (tarray tdouble k) VAL val_p;
      data_at_ Ews (Tstruct _csr_matrix noattr) q;
      data_at sh t_coo
       (rp, (cp, (vp,
        (Vint (Int.repr (Zlength (coo_entries coo'))),
         (Vint (Int.repr maxn),
          (Vint (Int.repr (coo_rows coo')), Vint (Int.repr (coo_cols coo')))))))) p;
       data_at sh (tarray tuint maxn)
         (map (fun e : $\mathbb{Z}$ * $\mathbb{Z}$ * ftype Tdouble => Vint (Int.repr (fst (fst e))))
           (coo_entries coo') ++ Zrepeat Vundef (maxn - Zlength (coo_entries coo')))
         rp;
       data_at sh (tarray tuint maxn)
         (map (fun e : $\mathbb{Z}$ * $\mathbb{Z}$ * ftype Tdouble => Vint (Int.repr (snd (fst e))))
           (coo_entries coo') ++ Zrepeat Vundef (maxn - Zlength (coo_entries coo')))
         cp;
       data_at sh (tarray tdouble maxn)
         (map (fun e : $\mathbb{Z}$ * $\mathbb{Z}$ * float => Vfloat (snd e)) (coo_entries coo') ++
          Zrepeat Vundef (maxn - Zlength (coo_entries coo'))) vp).
\end{lstlisting}

\section{Proving the partial\_CSR relation}
Now we need to prove that a relation exists that satisfies the lemmas
in \autoref{lst:axioms} and \autoref{lst:axioms2}.
\begin{lstlisting}
Derfinition partial_CSR {t: type} (i: $\mathbb{Z}$) (r: $\mathbb{Z}$) (coo: coo_matrixt t)
             (ROWPTR COLIND VAL: list C.val) : Prop :=
     (* fill in a definition here *).     
\end{lstlisting}
\pagebreak
In defining this relation, the concept of ``COO matrix up to entry $i$''
will be useful:

\begin{lstlisting}
Definition coo_upto (i: $\mathbb{Z}$) {t} (coo: coo_matrix t) : coo matrix t :=
 {| coo_rows := coo_rows coo;
    coo_cols := coo_cols coo;
    coo_entries := sublist 0 i (coo_entries coo)
 |}.
\end{lstlisting}
We say a \lstinline{coo_matrix} is well-formed if,
in every coordinate-tuple $(r_i,c_i,x)$
we have $0 \le r_i < \mathrm{coo\_rows}$ and
$0 \le c_i < \mathrm{coo\_cols}$.  Clearly, if \lstinline{coo} is well-formed,
then \lstinline{coo_upto $i$ coo} is well-formed.

\autoref{fig:CSR-partial} suggests that a partial COO matrix
(such as \lstinline{coo_upto 4 coo}) should relate to a partial CSR matrix
(such as the one shown in the figure).  But a partial CSR matrix is
not well-formed, because the \lstinline{row_ptr} array is not filled in.
However, the partial CSR matrix in that figure can be trivially completed by
filling in all the empty slots with 3, i.e., the number of elements
in the \lstinline{val} array.

Therefore the heart of the \lstinline{partial_CSR} relation should
be a relation between a \emph{complete} COO matrix (namely,
 \lstinline{coo_upto i coo}) and a \emph{complete} CSR matrix.
 We call this relation \lstinline{coo_csr}:
\begin{lstlisting}
Definition coo_csr {t} (coo: coo_matrix t) (csr: csr_matrix t) : Prop := . . .
\end{lstlisting}

The most natural \emph{semantic}
relation is that both of these represent the same
matrix:
\begin{lstlisting}
(* wrong *)  exists m: matrix t, coo_to_matrix coo m $~$/\ csr_to_matrix csr m.
 \end{lstlisting}
However, it turns out that this relation is not strong enough to
support the induction. (At least, I don't have evidence that it
is strong enough.)  The reason is that either matrix might contain
explicit zero values, and the COO matrix might have an $(i,j,0)$
where the CSR matrix has no entry, or vice versa.  Therefore,
we need a more structural assurance that if the 
COO matrix has one or more entries at $(i,j)$ then the CSR
matrix has a corresponding entry, and vice versa.

The relation is therefore defined as,
\begin{lstlisting}
Inductive coo_csr {t} (coo: coo_matrix t) (csr: csr_matrix t) : Prop :=
 build_coo_csr: forall
    (coo_csr_rows: coo_rows coo = csr_rows csr)
    (coo_csr_cols: coo_cols coo = csr_cols csr)
    (coo_csr_vals: Zlength (csr_vals csr) = count_distinct (coo_entries coo))
    (coo_csr_entries: entries_correspond coo csr)
    (coo_csr_zeros: no_extra_zeros coo csr),
    coo_csr coo csr.
\end{lstlisting}
That is, the COO and CSR matrix must have the same number of columns;
the length of the CSR \lstinline{val} array must be the
number of distinct $(i,j)$ coordinates of the COO;
for every entry in the COO matrix there must be
a corresponding element of the CSR matrix (entries\_\linebreak[0]correspond);
and for every element in the CSR matrix there must be a corresponding
entry in the COO matrix (no\_\linebreak[0]extra\_\linebreak[0]zeros).
Those latter relations are defined as,
\begin{lstlisting}
Definition entries_correspond {t} (coo: coo_matrix t) (csr: csr_matrix t) :=
forall h, 
0 $\le$ h < Zlength (coo_entries coo) ->
let '(r,c) := fst (Znth h (coo_entries coo)) in
let k := cd_upto (h+1) coo - 1 in
  Znth r (csr_row_ptr csr) $\le$ k < Znth (r+1) (csr_row_ptr csr) /\
  Znth k (csr_col_ind csr) = c /\
  sum_any (map snd (filter (coord_eqb (r,c) oo fst) (coo_entries coo)))
          (Znth k (csr_vals csr)).

Definition no_extra_zeros {t} (coo: coo_matrix t) (csr: csr_matrix t) := 
  forall r k, 0 $\le$ r < coo_rows coo ->
     Znth r (csr_row_ptr csr) $\le$ k < Znth (r+1) (csr_row_ptr csr) ->
     let c := Znth k (csr_col_ind csr) in  In (r,c) (map fst (coo_entries coo)).
\end{lstlisting}

Having established this \lstinline{coo_csr}
relation of a complete COO matrix to a complete CSR matrix,
we can use it in relating partial matrices,
as shown in \autoref{lst:partial-csr-def}.

\begin{lstlisting}[language=Coq,float,basicstyle=\ttfamily\footnotesize,
captionpos=b,caption=Wellformedness of CSR matrices and definition of \lstinline{partial_CSR},label=lst:partial-csr-def]
Inductive csr_matrix_wellformed {t} (csr: csr_matrix t) : Prop :=
 build_csr_matrix_wellformed:
 forall (CSR_wf_rows: 0 $\le$ csr_rows csr)
    (CSR_wf_cols: 0 $\le$ csr_cols csr)
    (CSR_wf_vals: Zlength (csr_vals csr) = Zlength (csr_col_ind csr))
    (CSR_wf_vals': Zlength (csr_vals csr) = Znth (csr_rows csr) (csr_row_ptr csr))
    (CSR_wf_sorted: list_solver.sorted Z.le (0 :: csr_row_ptr csr ++ [Int.max_unsigned]))
    (CSR_wf_rowsorted: forall r, 0 $\le$ r < csr_rows csr -> 
         sorted Z.lt 
          (-1 :: sublist (Znth r (csr_row_ptr csr)) (Znth (r+1) (csr_row_ptr csr)) 
                         (csr_col_ind csr)
           ++ [csr_cols csr])),
 csr_matrix_wellformed csr.

Inductive partial_CSR (h: $\mathbb{Z}$) (r: $\mathbb{Z}$) (coo: coo_matrix Tdouble)
      (rowptr: list val) (colind: list val) (val: list val) : Prop :=
build_partial_CSR: forall 
  (partial_CSR_coo: coo_matrix_wellformed coo)
  (partial_CSR_coo_sorted: sorted coord_le (coo_entries coo))
  (partial_CSR_i: 0 $\le$ h $\le$ Zlength (coo_entries coo))
  (partial_CSR_r: -1 $\le$ r $\le$ coo_rows coo)
  (partial_CSR_r': Forall (fun e => fst (fst e) $\le$ r) (coo_entries (coo_upto h coo)))
  (partial_CSR_r'': Forall (fun e => fst (fst e) $\ge$ r)
                 (sublist h (Zlength (coo_entries coo)) (coo_entries coo)))
  (csr: csr_matrix Tdouble)
  (partial_CSR_wf: csr_matrix_wellformed csr)
  (partial_CSR_coo_csr: coo_csr (coo_upto h coo) csr)
  (partial_CSR_val: sublist 0 (Zlength (csr_vals csr)) val = map Vfloat (csr_vals csr))
  (partial_CSR_colind: sublist 0 (Zlength (csr_col_ind csr)) colind =
                       map (Vint oo Int.repr) (csr_col_ind csr))
  (partial_CSR_rowptr: sublist 0 (r+1) rowptr =
                       map (Vint oo Int.repr) (sublist 0 (r+1) (csr_row_ptr csr)))
  (partial_CSR_val': Zlength val = count_distinct (coo_entries coo))
  (partial_CSR_colind': Zlength colind = count_distinct (coo_entries coo))
  (partial_CSR_rowptr': Zlength rowptr = coo_rows coo + 1)
  (partial_CSR_dbound: count_distinct (coo_entries coo) $\le$ Int.max_unsigned),
  partial_CSR h r coo rowptr colind val.
\end{lstlisting}

With this definition, proofs of all the \lstinline{partial_CSR}
lemmas proceed in a straightforward (though tedious) manner.

\section{Conclusion}

The program is not long, but its invariants are suprisingly intricate.
The definitions, properties, and lemmas for
\lstinline{coo_csr} and 
\lstinline{partial_CSR}
took 1571 lines of Coq, and based on those the VST proof took 412 lines.
In the process of developing the proof I found and fixed five bugs
in the program---four off-by-one errors and one in tricky situation
(discussed above) 
resulting from initializing the unsigned integer variable \lstinline{r}
to $-1$ (modulo 2 to the wordsize).

This proved-correct module is now ready to serve its intended purpose:
as a composable verified component in any verified program---such
as finite-element PDE solution---that requires construction of Compressed
Sparse Row matrices.

\vspace{\baselineskip}
\noindent\textbf{Acknowledgments.}  
This research was supported in part by NSF Grant CCF-2219757.  I thank David Bindel
for explanations and advice regarding construction and applications of sparse matrices.

\bibliographystyle{eptcs}
\bibliography{appel}

\begin{thebibliography}{1}
\providecommand{\bibitemdeclare}[2]{}
\providecommand{\surnamestart}{}
\providecommand{\surnameend}{}
\providecommand{\urlprefix}{Available at }
\providecommand{\url}[1]{\texttt{#1}}
\providecommand{\href}[2]{\texttt{#2}}
\providecommand{\urlalt}[2]{\href{#1}{#2}}
\providecommand{\doi}[1]{doi:\urlalt{https://doi.org/#1}{#1}}
\providecommand{\eprint}[1]{arXiv:\urlalt{https://arxiv.org/abs/#1}{#1}}
\providecommand{\bibinfo}[2]{#2}

\bibitemdeclare{inproceedings}{appel24:vcfloat}
\bibitem{appel24:vcfloat}
\bibinfo{author}{Andrew \surnamestart Appel\surnameend} \&
  \bibinfo{author}{Ariel \surnamestart Kellison\surnameend}
  (\bibinfo{year}{2024}): \emph{\bibinfo{title}{VCFloat2: Floating-Point Error
  Analysis in Coq}}.
\newblock In: {\slshape \bibinfo{booktitle}{Proceedings of the 13th ACM SIGPLAN
  International Conference on Certified Programs and Proofs}},
  \bibinfo{series}{CPP 2024}, \bibinfo{publisher}{Association for Computing
  Machinery}, \bibinfo{address}{New York, NY, USA}, p.
  \bibinfo{pages}{14–29}, \doi{10.1145/3636501.3636953}.

\bibitemdeclare{unpublished}{VerifiableCmanual}
\bibitem{VerifiableCmanual}
\bibinfo{author}{Andrew~W. \surnamestart Appel\surnameend},
  \bibinfo{author}{Lennart \surnamestart Beringer\surnameend},
  \bibinfo{author}{Qinxhiang \surnamestart Cao\surnameend} \&
  \bibinfo{author}{Josiah \surnamestart Dodds\surnameend}
  (\bibinfo{year}{2019}): \emph{\bibinfo{title}{{V}erifiable {C}: applying the
  {V}erified {S}oftware {T}oolchain to {C} programs}}.
\newblock \bibinfo{note}{{\tt https://vst.cs.princeton.edu/download/VC.pdf}}.

\bibitemdeclare{book}{barret94:templates}
\bibitem{barret94:templates}
\bibinfo{author}{Richard \surnamestart Barrett\surnameend},
  \bibinfo{author}{Michael \surnamestart Berry\surnameend},
  \bibinfo{author}{Tony~F. \surnamestart Chan\surnameend},
  \bibinfo{author}{James \surnamestart Demmel\surnameend},
  \bibinfo{author}{June \surnamestart Donato\surnameend}, \bibinfo{author}{Jack
  \surnamestart Dongarra\surnameend}, \bibinfo{author}{Victor \surnamestart
  Eijkhout\surnameend}, \bibinfo{author}{Roldan \surnamestart Pozo\surnameend},
  \bibinfo{author}{Charles \surnamestart Romine\surnameend} \&
  \bibinfo{author}{Henk \surnamestart van~der Vorst\surnameend}
  (\bibinfo{year}{1994}): \emph{\bibinfo{title}{Templates for the Solution of
  Linear Systems: Building Blocks for Iterative Methods}}.
\newblock \bibinfo{publisher}{SIAM}.

\bibitemdeclare{article}{vst-floyd}
\bibitem{vst-floyd} \bibinfo{author}{Qinxiang \surnamestart
  Cao\surnameend}, \bibinfo{author}{Lennart \surnamestart
  Beringer\surnameend}, \bibinfo{author}{Samuel \surnamestart
  Gruetter\surnameend}, \bibinfo{author}{Josiah \surnamestart
  Dodds\surnameend} \& \bibinfo{author}{Andrew~W. \surnamestart
  Appel\surnameend} (\bibinfo{year}{2018}):
  \emph{\bibinfo{title}{{VST-Floyd}: A Separation Logic Tool to Verify
    Correctness of {C} Programs}}.  \newblock {\slshape
    \bibinfo{journal}{J. Autom. Reason.}}
  \bibinfo{volume}{61}(\bibinfo{number}{1-4}),
  pp. \bibinfo{pages}{367--422},
  \doi{10.1007/s10817-018-9457-5}.

\bibitemdeclare{inproceedings}{kellison23:laproof}
\bibitem{kellison23:laproof}
\bibinfo{author}{Ariel~E. \surnamestart Kellison\surnameend},
  \bibinfo{author}{Andrew~W. \surnamestart Appel\surnameend},
  \bibinfo{author}{Mohit \surnamestart Tekriwal\surnameend} \&
  \bibinfo{author}{David \surnamestart Bindel\surnameend}
  (\bibinfo{year}{2023}): \emph{\bibinfo{title}{LAProof: A Library of Formal
  Proofs of Accuracy and Correctness for Linear Algebra Programs}}.
\newblock In: {\slshape \bibinfo{booktitle}{2023 IEEE 30th Symposium on
  Computer Arithmetic (ARITH)}}, pp. \bibinfo{pages}{36--43},
  \doi{10.1109/ARITH58626.2023.00021}.

\end{thebibliography}

\end{document}